\documentclass{article}

\usepackage{amsfonts}
\usepackage{amsmath}
\usepackage{amssymb}
\usepackage{amsthm}
\usepackage{ifthen}

\newtheorem{theorem}{Theorem}[section]
\newtheorem{proposition}[theorem]{Proposition}
\newtheorem{corollary}[theorem]{Corollary}
\newtheorem{lemma}[theorem]{Lemma}
\theoremstyle{remark}
\newtheorem{remark}[theorem]{Remark}

\theoremstyle{definition}

\newcommand{\parder}[3][Default]{
	\frac{\partial \ifthenelse{\equal{#1}{Default}}{}{^{#1}}#2}{
              \partial #3 \ifthenelse{\equal{#1}{Default}}{}{^{#1}}}}

\newcommand{\jac}{\operatorname{\mathcal J}}

\newcommand{\diag}{\operatorname{diag}}
\newcommand{\GL}{\operatorname{GL}}

\newcommand{\C}{{\mathbb C}}
\newcommand{\F}{{\mathbb F}}
\newcommand{\R}{{\mathbb R}}

\newcommand{\Z}{{\mathbb Z}}
\newcommand{\tp}{^{\rm t}}

\newcommand{\rk}{\operatorname{rk}}
\newcommand{\Mat}{\operatorname{Mat}}
\newcommand{\chr}{\operatorname{chr}}

\title{Some remarks on the Jacobian conjecture and polynomial endomorphisms}
\author{Dan Yan \\
School of Mathematical Sciences, Graduate University of \\
Chinese Academy of Sciences, Beijing 100049, China \\
\emph{E-mail:} yan-dan-hi@163.com
\and
Michiel de Bondt\footnote{The second author was supported by the Netherlands
                          Organisation for Scientific Research (NWO).} \\
Department of Mathematics, Radboud University \\
Nijmegen, The Netherlands \\
\emph{E-mail:} M.deBondt@math.ru.nl}

\allowdisplaybreaks

\begin{document}

\maketitle

\begin{abstract}
\noindent
In this paper, we first show that homogeneous Keller maps are injective on lines
through the origin. We subsequently formulate a generalization, which is
that under some conditions, a polynomial endomorphism with $r$ homogeneous parts of
positive degree does not have $r$ times the same image point on a line through
the origin, in case its Jacobian determinant does not vanish anywhere on that
line. As a consequence, a Keller map of degree $r$ does not take the same values
on $r > 1$ collinear points, provided $r$ is a unit in the base field.

Next, we show that for invertible maps $x + H$ of degree $d$,
such that $\ker \jac H$ has $n-r$ independent vectors over the base field,
in particular for invertible power linear maps $x + (Ax)^{*d}$ with
$\rk A = r$, the degree of the inverse of $x + H$ is at most $d^r$.
\end{abstract}


\section{Introduction}

Throughout this paper, we will write $K$ for any field and
$K[x] = K[x_1,x_2,\ldots,\allowbreak x_n]$ for the polynomial algebra over $K$ with $n$
indeterminates $x = x_1, x_2, \ldots, x_n$. Let $F = (F_1, F_2, \ldots, F_m): K^n
\rightarrow K^m$ be a polynomial mapping, i.e.\@ $F_i \in K[x]$ for all $i \le m$.
Denote by $\jac F$ the $(m \times n)$-matrix
$$
\jac F = \left( \begin{array}{cccc}
\parder{}{x_1} F_1 & \parder{}{x_2} F_1 & \cdots & \parder{}{x_n} F_1 \\
\parder{}{x_1} F_2 & \parder{}{x_2} F_2 & \cdots & \parder{}{x_n} F_2 \\
\vdots & \vdots &  & \vdots \\
\parder{}{x_1} F_m & \parder{}{x_2} F_m & \cdots & \parder{}{x_n} F_m
\end{array} \right)
$$
The well-known Jacobian conjecture (JC), raised by O.H. Keller in 1939 in
\cite{keller}, states that in case the characteristic $\chr K$ of $K$ is
zero, a polynomial mapping $F: K^n \rightarrow K^n$ is invertible if the
Jacobian determinant $\det \jac F$ is a nonzero constant.\footnote{The actual
formulation by Keller, also known as Keller's conjecture, was over $\Z$ instead of
$K$, with $\det \jac F \in \{-1,1\}$ a unit in $\Z$, but Keller's conjecture and the
JC are equivalent when they are quantified over all dimensions.}
This conjecture has been attacked by many people from various research fields,
but is still open, even for $n = 2$. Only the case $n=1$ is obvious, but the map
$x_1 \mapsto x_1 - x_1^q$ over $\F_q$ (which is the zero map)
shows that $\chr K = 0$ is required. For more information
about the wonderful 70-year history, see \cite{bcw}, \cite{arnobook}, and the
references therein. For more recent developments, see the second author's
Ph.D. thesis \cite{homokema} and the survey article \cite{amazing}.

Among the vast interesting and valid results, a satisfactory positive result
was obtained by S.S.S. Wang in 1980 in \cite{wang}, which is that the Jacobian
conjecture holds when the degree of the concerned polynomial map is equal to two.
A more simple proof of Wang's result was obtained by showing that a quadratic Keller map is
injective over the algebraic closure of the base field $K \ni \frac12$, since it
had already been shown that that is sufficient for concluding invertibility, see
\cite[Prop.\@ 17.9.6]{grothendieck}, \cite[Lm.\@ 3]{yagzhev}, \cite{cynkrusek} and
\cite{rudin}. There are other authors that can be added to this list, but some authors
only proved bijectivity.

The proof in \cite[(2.4)]{bcw} and \cite[Prop.\@ 4.3.1]{arnobook} of Wang's injectivity
result is due to S. Oda and is roughly as follows. Assume the opposite, say that $F(a) = F(b)$.
The first step is reducing to the case that $b = 0$. The second step is showing that $a = 0$
in case $F(a) = F(0)$. The second step can easily be generalized in the sense that for Keller
maps with terms of degree $0$, $1$, and $d$ only, we have $a = 0$ in case $F(a) = F(0)$,
where $d \ge 2$ is arbitrary. But the next proposition shows a stronger statement
for the line $Ka$ through $a$ and the origin.

\begin{proposition} \label{genprop}
Assume $F: K^n \rightarrow K^n$ is a Keller map with terms of degree $0$, $1$, and $d$
only, where $d > 1$ is a unit in $K$, and $a \in K^n$. Then $F|_{Ka}$ is injective.
\end{proposition}

\begin{proof}
By replacing $K$ by its algebraic closure, we may assume that $K$ is algebraically
closed. Write $F = F^{(0)} + F^{(1)} + F^{(d)}$, where $F^{(i)} \in K[x]^n$ has terms of
degree $i$ only. By Euler's theorem for homogeneous functions, we have $F^{(d)}_i =
d^{-1}\sum_{j=1}^n x_j\parder{}{x_j} F^{(d)}_i = d^{-1} \jac F^{(d)}_i \cdot x$
for all $i$. Hence $F(a) = F(b)$, if and only if
$$
F^{(0)} + F^{(1)}(a) + d^{-1}(\jac F^{(d)})|_{a} \cdot a
= F^{(0)} + F^{(1)}(b) + d^{-1}(\jac F^{(d)})|_{b} \cdot b
$$
Assume $F(a) = F(\lambda a)$ for some $\lambda \in K$.
By subtracting $F^{(0)}$ and substituting $b = \lambda a$, we get
$$
\Big(\jac F^{(1)}|_{a} + d^{-1}(\jac F^{(d)})|_{a}\Big) \cdot a
= \Big(\lambda \jac F^{(1)}|_{a} + d^{-1}\lambda^d(\jac F^{(d)})|_{a}\Big) \cdot a
$$
where $|_f$ means substituting $f$ for $x$.
Since $\jac F^{(1)}|_{a} = \jac F^{(1)}$, this is equivalent to
$$
\Big((1-\lambda)\jac F^{(1)} + d^{-1}(1-\lambda^d)
(\jac F^{(d)})|_{a}\Big) \cdot a = 0
$$
Notice that $d^{-1}(1-\lambda^d) = d^{-1}(1 + \lambda + \cdots + \lambda^{d-1})
\cdot (1-\lambda)$, and $\jac F^{(d)}$ is homogeneous of degree $d-1$.
If we define
$\mu := \sqrt[d-1]{d^{-1}(1 + \lambda + \cdots + \lambda^{d-1})}$, then the above
is equivalent to
$$
(1-\lambda) \Big(\jac F^{(1)} + (\jac F^{(d)})|_{\mu a}\Big) \cdot a = 0
$$
i.e.\@ $(\jac F)|_{ \mu a} \cdot (\lambda-1)a = 0$. Since
$\det (\jac F)|_{ \mu a}$ is a nonzero constant, $(\lambda-1)a = 0$ is
the only possibility. Hence $b = a$, as desired.
\end{proof}

\noindent
The map $F = x_1 - x_1^q$ over $\F_q$ shows that the condition that
$d$ is a unit in $K$ is required in proposition \ref{genprop}.
We will generalize proposition \ref{genprop} in section \ref{gensec}.
In the above proof of \ref{genprop}, we do not need that $\det \jac F$ is
a nonzero constant: it is sufficient that $\det \jac F$ does not vanish on
$Ka$, provided we can take $(d-1)$-th roots in $K$. This is expressed in
theorem \ref{genth} in section \ref{gensec}.

Although it is sufficient to show injectivity over the algebraic closure,
it is not true that injective polynomial maps are automatically
invertible: take for instance the map $F = x_1 + x_1^3$ over $\R$.
Injective quadratic maps do not need to be invertible either: take
for instance the map $F = (x_1 + x_2 x_3, x_2 - x_1 x_3, x_3)$ over $\R$,
which is injective because for fixed $x_3$ it corresponds to the linear
map $\tilde{F} = (x_1 + x_2 x_3, x_2 - x_1 x_3)$ over $\R[x_3]$,
which has Jacobian determinant $1 + x_3^2 \ne 0$ just as $F$ itself.
It is however true that injective polynomial maps over $\R$ are automatically
surjective, see \cite{bbr}.

Quadratic polynomial maps over a field $K \ni \frac12$,
such that Jacobian determinant does not
vanish anywhere, are always injective. This follows from both Oda's proof of
Wang's theorem and the proof of proposition \ref{genprop}.
In the proof of proposition \ref{genprop}, the Keller condition is used to replace
$K$ by its algebraic closure without affecting that the Jacobian determinant
is nonzero everywhere, but since a $(d-1)$-th root is taken, there
is no need to replace $K$ by its algebraic closure when $d = 2$.
For twodimensional Keller maps over an algebraically closed field of characteristic
zero, there is a nice result due to J. Gwo{\'z}dziewicz in \cite{gw}, namely that
they are invertible  when they are injective on one single line (any line).
His short proof makes use of the Abhyankar-Moh-Suzuki theorem, see e.g.\@
\cite[Th.\@ 5.3.5]{arnobook} for the characteristic zero case which is used
in \cite{gw}.

However in \cite{pinchuk}, S. Pinchuk constructed a twodimensional
polynomial map over $\R$ of degree 25
with nonzero Jacobian everywhere, which is {\em not injective}\footnote{Pinchuk
claims to have constructed a map of degree 40, but does not give it explicitly.
Inspection of his proofs leads to a map of degree $\ge 25$ which can be reduced to
degree $= 25$ by way of an elementary automorphism.}
The counterexample of Pinchuk can be transformed to a polynomial map of
degree three in larger dimension by way of reduction techniques that
are indicated below. In that manner, in 1999 E. Hubbers constructed
a non-injective cubic polynomial map of a very special form
in dimension 1999 over $\R$ with nonzero Jacobian everywhere, namely a so-called
Dru{\.z}kowski map. This result has not been published, though.

After having solved the quadratic case of the JC, it is a natural
question what happens when the degree of the polynomial map is equal to three.
This case of the JC has not been solved yet, but instead,
it is proved that the JC holds in general in case it holds
for cubic maps of the form $x + H$, where $x = (x_1, x_2, \ldots, x_n)$ is the
identity map in dimension $n$, and $H = (H_1, H_2, \ldots, H_n)$ is homogeneous
of degree three, i.e.\@ each $H_i$ is either homogeneous of degree three or
zero.

In fact, the JC holds in general if for some $d \ge 3$ (any $d \ge 3$), the JC
holds for maps of the form $x + H$, where $x = (x_1, x_2, \ldots, x_n)$ is the
identity map in dimension $n$, and $H = (H_1, H_2, \ldots, H_n)$ is homogeneous
of degree $d$, and proposition \ref{genprop} shows that Keller maps of this
form are injective on lines through the origin.

A subsequent reduction is due to L.M. Dru{\.z}kowski, which asserts that the
JC is true in general if it is true for maps of the form $x + H$ such that
$H_i$ is either a third power of a linear form or zero for each $i$.
Therefore, a polynomial map $x + H$ as such is called a {\em Dru{\.z}kowski map}.
If $A \in \Mat_n(K)$, then we can take the matrix product $Ax$ of $A$
with the column vector $x$, and $Ax$ is a column vector of linear forms.
Next, we can take the Hadamard product $(Ax) * (Ax) * (Ax)$, which is called
the third Hadamard power of $Ax$ and denoted as $(Ax)^{*3}$ by many
authors\footnote{In fact, it is popular to denote $\alpha^{\diamond d}$
for the `composition' of $d$ copies of some object $\alpha$ by means of applying
any commutative binary operator $\diamond$ exactly $d-1$ times, such as $K^{\times n}$
for the $n$-dimensional vector space over $K$, and in an invertible context,
the JC is about the existence of $F^{\circ(-1)}$.}. Thus
Dru{\.z}kowski showed that the JC is true in general if it is true for maps
of the form $x + (Ax)^{*3}$.

Similarly, one can define $x + (Ax)^{*d}$, and just as above,
the JC is true in general in case it is true for maps of the
form $x + (Ax)^{*d}$, where $d$ is one's favorite integer larger than two.
If $F$ is an invertible polynomial map of degree $d$ in dimension $n$, then the
degree of its inverse is at most $d^{n-1}$. This has been proved
in \cite[Cor.\@ (1.4)]{bcw}, which is a direct consequence of a more or less
similar result about birational maps in projective space by O. Gabber,
see \cite[Th.\@ (1.5)]{bcw}.
But if $F$ is of the form $x + (Ax)^{*d}$ and $\rk A = r$, then
the degree of the inverse of $F$ is at most $d^r$. We will prove this
in section \ref{druzsec}.

\section{Polynomial maps that take the same values on several collinear points}
\label{gensec}

At first, a lemma with a generalized Vandermonde matrix which is assumed to have
full rank. The powers in the Vandermonde matrix correspond to the degrees
of the homogeneous parts of the polynomial map in theorem \ref{genth},
the main theorem of this section.

\begin{lemma} \label{genlm}
Let $K$ be a field and assume that $G : K \rightarrow K^m$ is given by
$$
G(t) = C (t^{d_1}, t^{d_2}, \ldots, t^{d_r}, t^{d_{r+1}})
$$
where $C \in \Mat_{m,r+1}(K)$. Suppose that $G(a_1) = G(a_2) = \cdots =
G(a_r) = 0$ and that $\rk A^{(r,r)} = r$, where $A^{(s,r)}$ is
the generalized Vandermonde matrix
\begin{equation} \label{genVan}
A^{(s,r)} := \left( \begin{array}{cccc}
a_1^{d_1} & a_2^{d_1} & \cdots & a_r^{d_1} \\
a_1^{d_2} & a_2^{d_2} & \cdots & a_r^{d_2} \\
\vdots & \vdots & \ddots & \vdots \\
a_1^{d_s} & a_2^{d_s} & \cdots & a_r^{d_s}
\end{array} \right)
\end{equation}
Then $\rk C \le 1$ and the last column of $C$ is nonzero in case $C \ne 0$.
\end{lemma}

\begin{proof}
Since $\rk A^{(r+1,r)} = r$, there is only one relation between the rows of $A$
up to scalar multiplication, say $v\tp A^{(r+1,r)} = 0$ for some nonzero vector
$v \in K^{r+1}$. By assumption, we have $C A^{(r+1,r)} = 0$, whence each row of
$C$ is a scalar multiple of $v\tp$. Thus $\rk C \le 1$.

If $C \ne 0$ and the last column of $C$
is zero, then $C A^{(r+1,r)} = 0$ would imply $\rk A^{(r,r)} < r$, which is a
contradiction.
\end{proof}

\begin{theorem} \label{genth}
Let $K$ be a field and assume $F: K^n \rightarrow K^m$ is a polynomial
map such that the degree of each term of $F$ is contained in $\{d_1, d_2, \ldots,
d_r, d_{r+1}\} \ni 0$. Assume furthermore that for all $c \in K^r$,
$$
\parder{}{t} (t^{d_{r+1}} + c_r t^{d_r} + \cdots + c_2 t^{d_2} + c_1 t^{d_1})
$$
has a root in $K$.

Take $b \in K^n$ nonzero and assume $F(a_1 b) = F(a_2 b) = \cdots = F(a_r b)$
and $\rk A^{(r,r)} = r$ for certain $a_i \in K$, where $A^{(r,r)}$ is
as in (\ref{genVan}). Then there is an $a_{r+1} \in K$ such that
$(\jac F)|_{a_{r+1}b} \cdot b = 0$. In particular,
$\rk (\jac F)|_{a_{r+1}b} < n$.
\end{theorem}

\begin{proof}
In this proof, we will write $|_{x = f}$ and $|_{t = f}$ for the substitution of $f$ for $x$ and $t$
respectively. Take $G(t) = F(tb) - F(a_1 b)$. Then $G(a_i) = F(a_i b) - F(a_1 b) = 0$ for
each $i \le r$, and the degree of each term of $G(t)$ is contained in $\{d_1, d_2, \ldots,
d_r, d_{r+1}\}$. Hence there exists a $C$ as in lemma \ref{genlm}, and by the chain rule
\begin{equation} \label{FGC}
(\jac F)|_{x = tb} \cdot b = \jac_t F(tb) = \jac_t G(t) =
C \jac_t(t^{d_1}, t^{d_2}, \ldots, t^{d_r}, t^{d_{r+1}})
\end{equation}
If $\jac_t G(t) = 0$, then we can take $a_{r+1}$ arbitrary. Hence
assume that $\jac_t G(t) \ne 0$, say that
$$
\parder{}{t} G_i(t) = C_i \jac_t(t^{d_1}, t^{d_2}, \ldots, t^{d_r}, t^{d_{r+1}}) \ne 0
$$
Since $\rk C = 1$ and the last column of $C$ is nonzero, every column of $C$ is dependent
of the last one. Therefore, $C_i \ne 0$ refines to $C_{i(r+1)} \ne 0$.
By assumption, $C_{i(r+1)}^{-1} \parder{}{t} G_i(t)$ has a root $a_{r+1} \in K$. Thus
$$
C_i \big(\jac_t(t^{d_1}, t^{d_2}, \ldots, t^{d_r}, t^{d_{r+1}})\big)\big|_{t = a_{r+1}} =
C_{i(r+1)} \Big(C_{i(r+1)}^{-1} \parder{}{t} G_i(t)\Big)\Big|_{t=a_{r+1}} = 0
$$
Using that $C_i \ne 0$ and again that $\rk C = 1$, we see that every row of $C$ is
dependent of $C_i$, which gives $C (\jac_t(t^{d_1}, t^{d_2}, \ldots,
t^{d_r}, t^{d_{r+1}}))|_{t = a_{r+1}} = 0$. Substituting $t = a_{r+1}$ in (\ref{FGC})
subsequently gives $(\jac F)|_{x = a_{r+1}b} \cdot b = 0$, as desired.
\end{proof}

\begin{corollary} \label{gencr}
Let $K$ be a field and assume $F: K^n \rightarrow K^n$ is a polynomial
map such that the degree of each term of $F$
is contained in $\{d_1, d_2, \ldots, d_r, d_{r+1}\} \ni 0$, where
$\chr K \nmid d_{r+1} \ne 1$.

Assume $b \in K^n$ is nonzero, such that
$F(a_1 b) = F(a_2 b) = \cdots = F(a_r b)$ and $\rk A^{(r,r)} = r$
for certain $a_i \in K$, where $A^{(r,r)}$ is
as in (\ref{genVan}). Then $\det \jac F \notin K^{*}$.
\end{corollary}

\begin{proof}
By replacing $K$ by its algebraic closure, we may assume that $K$
is algebraically closed. The condition $\chr K \nmid d_{r+1} \ne 1$
tells us that for all $c \in K^r$,
$$
\parder{}{t} (t^{d_{r+1}} + c_r t^{d_r} + \cdots + c_2 t^{d_2} + c_1 t^{d_1})
$$
has a root in $K$. Hence by theorem \ref{genth}, there exists
an $a_{r+1} \in K$ such that  $(\jac F)|_{a_{r+1}b}$ does not have full rank.
Consequently, $x = a_{r+1}b$ is a root of $\det \jac F$. This gives the desired result.
\end{proof}

\noindent
If we take $r = 2$ in corollary \ref{gencr}, then
$d_3 > 1$ must be a unit in $K$ and $0 \in \{d_1,d_2\}$.
Furthermore, the conclusion of corollary \ref{gencr} is trivial when
$1 \notin \{d_1,d_2\}$, thus we may assume that $d_1 = 0$ and $d_2 = 1$.
Since the Vandermonde matrix $A^{(r,r)}$ always has full rank when $r = 2$
and $a_1 \ne a_2$, proposition \ref{genprop} is the case $r = 2$ of corollary
\ref{gencr}.

In corollary \ref{gencr2} below, we show that a Keller map of degree $r$
does not take the same values on $r$ collinear points, provided $r > 1$ is a
unit in the base field.

\begin{corollary} \label{gencr2}
Let $K$ be a field and assume $F: K^n \rightarrow K^n$ is a polynomial
map such that $F(p_1) = F(p_2) = \cdots = F(p_r)$ for distinct collinear
$p_i \in K^n$. If $\chr K \nmid r \ne 1$ and $r \ge \deg F$, then
$\det \jac F \notin K^{*}$.
\end{corollary}

\begin{proof}
By replacing $F$ by $F(x - p_1)$, we may assume that $p_i = a_i b$ for distinct
$a_i \in K$, where $b = p_2 - p_1 \ne 0$. Take $d_i = i-1$ for all $i \le r+1$. Then
$A^{(r,r)}$ is a regular Vandermonde matrix and
$\rk A^{(r,r)} = r$, where $A^{(r,r)}$ is as in (\ref{genVan}). Hence the desired
result follows from corollary \ref{gencr}.
\end{proof}

\begin{remark}
Just as with proposition \ref{genprop}, the map $f = x_1 - x_1^q$ with $r = 2$ and
$(d_1,d_2,d_3) = (0,1,q)$ shows that the condition that $\chr K \nmid d_{r+1}$ is
required in corollary \ref{gencr}. When we avoid the use of the condition that
$\chr K \nmid d_{r+1}$ by taking $r = 3$ and $(d_1,d_2,d_3,d_4) = (0,1,q,q+1)$,
the condition that the Vandermonde matrix $A^{(r,r)}$ has full rank appears necessary
in corollary \ref{gencr}.

The condition $\chr K \nmid d_{r+1}$ is also required if we take $r \ge \deg F$
and $d_i = i - 1$ for all $i$
in corollary \ref{gencr}, which we did in the proof of corollary \ref{gencr2}.
This is why the map $f = x_1 - x_1^q$ shows that the condition
$\chr K \nmid r \ne 1$ is required in corollary \ref{gencr2} as well.
\end{remark}

\section{The degree of the inverse of a polynomial map}
\label{druzsec}

We first show that the degree bound $d^{n-1}$ on the inverse of a polynomial map $F$ of degree $d$
in dimension $n$ is reached when $F = x + (Ax)^{*d}$ is power linear such that $\rk A = n-1$.

\begin{proposition}
Assume $F = x + (Ax)^{*d}$ is power linear of degree $d \ge 1$
such that $\rk A = n - 1$. If $\det \jac F = 1$, then
 $F$ is invertible and its inverse has degree $d^{n-1}$.
\end{proposition}

\begin{proof}
The case $d = 1$ is trivial, so assume that $d \ge 2$.
In \cite[Th.\@ 1]{chengquad} and \cite[Th.\@ 2.1]{druz}, it is shown that $F$ is linearly triangularizable
in case $d = 2$ and tame in case $d = 3$ respectively, but inspection of the proofs tells us that
in both cases, there exists a $T \in \GL_n(K)$ such that $T^{-1} F(Tx) = x + (Bx)^{*d}$, where
$B$ is lower triangular with zeroes on the diagonal. Furthermore, the case $d \ge 4$ follows in a
similar manner as the case $d = 3$, see also \cite[Th.\@ 4]{tong}.
Therefore, we may assume that $A$ is lower triangular with zeroes on the diagonal.

Let $G$ be the inverse of $F$. Then $G_1 = x_1$, and substituting $x = G$ in
$x_i = F_i - (A_{i1} x_1 + \cdots + A_{i(i-1)} x_{i-1})^d$ gives
$G_i = x_i - (A_{i1} G_1 + \cdots + A_{i(i-1)} G_{i-1})^d$, which is an inductive
formula for $G$.
Since $\rk A = n - 1$, we see that the entries $A_{i(i-1)}$ on the subdiagonal of $A$ are all
nonzero, and $\deg G_i = d \deg G_{i-1}$ follows for all $i$ by induction.
Hence $\deg G_n = d^{n-1} \deg G_1 = d^{n-1}$.
\end{proof}

\noindent
Before we prove our theorem, we first show that by way of a linear conjugation, we can
reduce to the case that $\ker \jac (F - x) \cap K^n = \{0\}^r \times K^{n-r}$.

\begin{lemma} \label{deglm}
Let $K$ be a field and assume $F: K^n \rightarrow K^n$ is an invertible polynomial map.
If $\ker \jac (F - x) \cap K^n$ has dimension $n-r$ as a $K$-space, then there exists
a $T \in \GL_n(K)$ such that for $G := T^{-1} F(Tx)$, we have
$$
\ker \jac (G - x) \cap K^n = \{0\}^r \times K^{n-r}
$$
Furthermore, $G$ is invertible and the degree of its inverse is the same as that of $F$.
\end{lemma}

\begin{proof}
Take $T \in \GL_n(K)$ such that the last $n-r$ columns of $T$ are a basis of the
$K$-space $\ker \jac (F - x) \cap K^n$, and set $G := T^{-1} F(Tx)$.
Then
$$
G - x = T^{-1} F(Tx) - T^{-1} Tx = T^{-1} (F - x)|_{Tx}
$$
and by the chain rule, $\jac (G - x) = T^{-1} \cdot
(\jac(F - x))|_{Tx} \cdot T$. Hence
\begin{align*}
\ker \jac (G - x) \cap K^n
&= T^{-1} \big(\!\ker \big(T^{-1} \cdot (\jac (F - x))|_{Tx}\big)\big) \cap K^n \\
&= T^{-1} (\ker \jac (F - x))|_{Tx} \cap T^{-1} (K^n)|_{Tx} \\
&= T^{-1} (\ker \jac (F - x) \cap K^n)|_{Tx} \\
&= T^{-1} (K T e_{r+1} + K T e_{r+2} + \cdots + K T e_{n})|_{\ldots} \\
&= K e_{r+1} + K e_{r+2} + \cdots + K e_{n} = \{0\}^r \times K^{n-r}
\end{align*}
where $e_i$ is the $i$-th standard basis unit vector.
If $\tilde{F}$ is the inverse of $F$, then $\tilde{G} := T^{-1} \tilde{F} (Tx)$ is
the inverse of $G$ and $\tilde{G}$ has the same degree as $\tilde{F}$, as desired.
\end{proof}

\begin{remark}
Write $\tilde{x} = (x_1,x_2,\ldots,x_r)$.
If we take for $B$ the first $r$ rows of $T^{-1}$ and for $C$ the first $r$ columns of
$T$, then we get $(G_1, G_2, \ldots, G_r) = B(F(Tx))$, and since $G_i \in K[\tilde{x}]$ for
all $i \le r$, the substitution $x_{r+1} = \cdots = x_n = 0$ has no effect, and
$(G_1, G_2, \ldots, G_r) = B(F(C\tilde{x}))$ follows. Furthermore, $BC$ is a leading principal
minor matrix of $T^{-1} T$ and hence equal to $I_r$, and since $T^{-1} \ker B =
\ker B T = \{0\}^r \times K^{n-r} = T^{-1} (\ker \jac (F - x) \cap K^n)|_{\ldots}$,
we have $\ker B = \ker \jac (F - x) \cap K^n$. Consequently, $(G_1, G_2, \ldots, G_r)$
and $F$ are paired in the sense of \cite{gz} when $r < n$, where the condition
$\ker B = \ker A$ is replaced by $\ker B = \ker \jac (F - x) \cap K^n$, to allow maps $F$
that are not of the form $x + (Ax)^{*d}$ as well.
\end{remark}

\begin{theorem} \label{degth}
Let $K$ be a field with $\chr K = 0$ and assume $F = x + H: K^n \rightarrow K^n$ is an
invertible polynomial map of degree $d$. If $\ker \jac H \cap K^n$
has dimension $n-r$ as a $K$-space, then the inverse polynomial map of $F$ has
degree at most $d^r$.
\end{theorem}

\begin{proof}
From lemma \ref{deglm}, it follows that by replacing $F$ by $T^{-1} F(Tx)$ for
a suitable $T \in \GL_n(K)$, we may assume that $\ker \jac H \cap K^n =
\{0\}^r \times K^{n-r}$. This means that
the last $r$ columns of $\jac H$ are zero. Hence $H_i \in K[x_1,x_2,\ldots,x_r]$
for all $i$. Let $x - G$ be the inverse of $F$ at this stage. Then
$G(F) = H$, because $x = (x - G)|_{F} = x + H - G(F)$.

Let $L = K(x_1,x_2,\ldots,x_r)$. Since $F_i, H_i \in L$ for all $i \le r$ and
$F_{r+1} ,F_{r+2}, \ldots, \allowbreak F_n$ are algebraically independent over $L$,
we obtain from $G_i(F) = H_i$ that $G_i \in K[x_1,x_2,\ldots,x_r]$ for all $i \le r$.
Hence $(x_1-G_1,x_2-G_2,\ldots,x_r-G_r)$ is the inverse polynomial map of
$(F_1,F_2,\ldots,F_r)$. Therefore, $\deg (x_1-G_1,x_2-G_2,\allowbreak \ldots,
x_r-G_r) \le d^{r-1}$ on account of \cite[Cor.\@ (1.4)]{bcw}.

Using that $F$ is the inverse of $x - G$, we obtain
by substituting $x = x - G$ in $G_i(F) = H_i$ that $G_i = H_i(x-G)$ for all $i$.
Since $\deg H_i \le d$ and $H_i \in K[x_1,\allowbreak x_2,\allowbreak \ldots,x_r]$
for all $i$, we get
$$
\deg G_i \le \deg H_i \cdot \deg (x_1-G_1,x_2-G_2,\ldots,x_r-G_r) \le d \cdot d^{r-1}
$$
for all $i$, which gives the desired result.
\end{proof}

\begin{corollary} \label{degcr}
Let $K$ be a field with $\chr K = 0$ and assume $F: K^n \rightarrow K^n$ is an
invertible polynomial map of the form $F = x + H$, where $H$
is power linear of degree $d$. If $\rk \jac H = r$,
then the inverse polynomial map of $F$ has degree at most $d^r$.
\end{corollary}

\begin{proof}
Say that $H = (Ax)^{*d}$. Then we have $\jac H = d \diag((Ax)^{*(d-1)}) \cdot A$,
where $\diag(v)$ stands for a square matrix with diagonal $v$ and
zeroes elsewhere. Consequently,
$\ker \jac H = \ker A$ and $\rk A = \rk \jac H = r$. Since
$A$ is a matrix over $K$, the dimension of $\ker A \cap K^n$ as a $K$-space is
$n-r$. Hence the desired result follows from theorem \ref{degth}.
\end{proof}

\noindent
Now that we have a better estimate of the degree of the inverse of a polynomial map,
it is interesting to know some cases that the polynomial map has an inverse on account
of the Keller condition. The following result was
proved by A. van den Essen in \cite[Prop.\@ 2.9]{arnobook} and C. Cheng in \cite[Th.\@ 2]{cheng}
for homogeneous maps and power linear maps.

\begin{theorem} \label{ec}
Let $K$ be a field with $\chr K = 0$ and write $\tilde{x} = (x_1,x_2,\ldots, x_r)$.
Assume that Keller maps $\tilde{F} = \tilde{x} + \tilde{H}$
in dimension $r$ over $K$ such that $\tilde{H}$ is (homogeneous) of degree $d$ are invertible.
Then Keller maps $F = x + H$ in dimension $n$ over $K$ such that $H$ is (homogeneous) of degree $d$
and $\ker \jac H \cap K^n$ has dimension $n - r$ are invertible as well.

A similar result holds when we replace `invertible' by `linearly triangularizable'.
\end{theorem}

\begin{proof}
Let $F$ be as above.
Just as in the proof of theorem \ref{degth}, we may assume that $H_i \in K[\tilde{x}]$ for all $i$.
Since $F$ is a Keller map in addition, the leading principal minor determinant of size $r$ of $\jac F$ is a
nonzero constant.
Now the leading principal minor matrix of size $r$ of $\jac F$ is of the form $\jac_{\tilde{x}}
\tilde{F}$ with $\tilde{F} = (F_1,F_2, \ldots, F_r)$, whence $\tilde{F}$
is a Keller map. By assumption, the inverse $\tilde{x} - \tilde{G}$
of $\tilde{F}$ exists, and $\tilde{F} - \tilde{G}(\tilde{F}) = \tilde{x}$.
If we define $G_i = \tilde{G}_i$ for all $i \le r$ and
$G_i = H_i(\tilde{x} - \tilde{G})$ for all $i \ge r+1$, then substituting $x = F$ in
$x_i - G_i$ gives
$$
F_i - G_i(F) = \left\{ \begin{array}{ll}
               \tilde{F}_i - \tilde{G}_i(\tilde{F}) = x_i & \mbox{if $i \le r$} \\
               x_i + H_i - H_i(\tilde{F} - \tilde{G}(\tilde{F})) = x_i & \mbox{if $i \ge r+1$}
\end{array} \right.
$$
Hence $x - G$ is the inverse of $F$.

In case $\tilde{F}$ is linearly triangularizable,
say that $\tilde{T}^{-1} \tilde{F}(\tilde{T}\tilde{x})$ has a lower triangular Jacobian,
then $T^{-1} F(Tx)$ has a lower triangular Jacobian as well if we define
$$
T = \left( \begin{array}{cc} \tilde{T} & \emptyset \\
    \emptyset & I_{n-r} \end{array} \right)
$$
which completes the proof of theorem \ref{ec}.
\end{proof}

\begin{corollary}
Let $K$ be a field with $\chr K = 0$ and assume that $F = x + H$ is a Keller map
in dimension $n$ over $K$ such that $\ker \jac H \cap K^n$ has
dimension $n - r$. Then $F$ is linearly triangularizable in the following cases:
\begin{enumerate}

\item[i)] $r \le 3$ and $H$ is homogeneous of degree $d \ge 2$,

\item[ii)] $r = 4$ and $H$ is homogeneous of degree $2$.

\end{enumerate}
Furthermore, $F$ is invertible as well in the following cases:
\begin{enumerate}

\item[iii)] $r = 2$ and $\deg H \le 101$,

\item[iv)] $r = 3$ and $\deg H \le 3$,

\item[v)] $r = 4$ and $H$ is homogeneous of degree $3$.

\end{enumerate}
\end{corollary}

\begin{proof}
This follows from theorem \ref{ec} and what is known about the Jacobian conjecture in small
dimensions $r$. For i), see \cite{dbvde}. For ii), use \cite{mo} and \cite[Th.\@ 7.4.4]{arnobook}.
For the rest, see \cite{arnobook} and the references therein.
\end{proof}

\end{document}